\documentclass[11pt]{article}

\usepackage[a4paper,margin=1in]{geometry}
\usepackage{amsmath,amssymb,amsthm,mathtools,bm,mathrsfs}
\usepackage{enumitem}
\usepackage{booktabs}
\usepackage{microtype}
\usepackage{indentfirst}
\usepackage{xcolor}
\usepackage{authblk}
\usepackage[colorlinks=true,linkcolor=blue!45!black,citecolor=blue!45!black,urlcolor=blue!45!black]{hyperref}
\usepackage[nameinlink,noabbrev]{cleveref}

\newtheorem{theorem}{Theorem}[section]
\newtheorem{lemma}[theorem]{Lemma}
\newtheorem{proposition}[theorem]{Proposition}

\theoremstyle{definition}

\newtheorem{remark}[theorem]{Remark}
\newtheorem{example}[theorem]{Example}

\numberwithin{equation}{section}

\newcommand{\Sph}{\mathbb S}
\newcommand{\Om}{\Omega}
\newcommand{\C}{\mathcal C}
\newcommand{\sig}{\sigma}
\newcommand{\dd}{\,\mathrm d}

\newcommand{\Id}{\mathrm{Id}}

\newcommand{\Hh}{\mathcal H}
\newcommand{\U}{\mathcal U}

\title{\bfseries The Capillary Dual Orlicz--Christoffel--Minkowski Problem}

\author[1]{Dongjing Ma}
\author[2]{Qiqi Zhang}
\author[3]{Hua Zhu}

\affil[1]{School of Mathematical Sciences, Qufu Normal University\\
Qufu, Shandong, 273165, People's Republic of China\\
\textit{Email: 3261019581@qq.com}}

\affil[2]{School of Mathematics and Information Science, Nanchang Hangkong University\\Nanchang,Jiangxi,330063, People's Republic of China\\
\textit{Email:71131@nchu.edu.cn}}

\affil[3]{School of Mathematics and Physics, Southwest University of Science and Technology\\
Mianyang, Sichuan, 621010, People's Republic of China\\
\textit{Email: zhuhmaths@mail.ustc.edu.cn}}

\date{}

\begin{document}

\maketitle
\vspace{-1.2em}

\begin{abstract}
We study a capillary dual Orlicz--Christoffel--Minkowski problem on a spherical cap.
For \(2\le k\le n\), \(1<p<k+1\), and \(0<\theta<\pi/2\), we consider
\[
 \sigma_k(\nabla^2h+h\sig)
 =
 f(\xi)h^{p-1}\frac{\rho^{k+1}}{\varphi(\rho)},
 \qquad
 \rho=(h^2+|\nabla h|^2)^{1/2},
\]
with the capillary Robin condition \(\nabla_\mu h=\cot\theta\,h\).
Assuming that \(f\) is positive, smooth and horizontally even,
\(g=f^{-1/(p+k-1)}\) satisfies the capillary spherical-convexity and boundary compatibility conditions, and
\(H(r)=\varphi(r)/r^p\) is strictly increasing and proper with
\(Q_\varphi(r)\le k+1\), we prove the existence of a smooth, positive, horizontally even solution that is strictly capillary convex.
The proof uses an Orlicz homotopy, uniform a priori estimates, a viscosity constant-rank argument when \(k<n\),
global oblique regularity, and degree theory.
\end{abstract}



\section{Introduction}

The Minkowski problem and the Christoffel--Minkowski problem are among the central problems in convex geometry.
Their \(L_p\), dual and Orlicz extensions have led to several related Brunn--Minkowski theories.
For general background on convex geometry and geometric reconstruction, see
\cite{Gardner1995,Schneider2014}.
Classical developments of the \(L_p\)-Minkowski problem can be found in
\cite{ChouWang2006,Lutwak1993,LYZ2004}, while the Christoffel--Minkowski problem for
Firey's \(p\)-sum was studied by Hu, Ma and Shen \cite{HuMaShen2004}.
Orlicz additions and the associated Brunn--Minkowski inequalities provide an extension of the power-based framework;
see \cite{GardnerEtAl2014,XiJinLeng2014}.

The dual theory and its geometric measures were developed in
\cite{HLYZ2016,Lutwak1975,Lutwak1988}.
The dual Orlicz theory provides a nonhomogeneous extension in which radial data are encoded through a general
positive function rather than a fixed power \cite{GardnerEtAl2015,Zhao2015,ZhuZhouXu2014}.
The dual Orlicz curvature measure and the corresponding dual Orlicz--Minkowski problem were introduced and studied
by Zhu, Xing and Ye \cite{ZXY2018}; related dual Minkowski and Orlicz--Minkowski problems can be found in \cite{HaberlEtAl2010,HuangHe2012,HLYZ2016,Li2014,Zhao2017}, while dual Orlicz Christoffel--Minkowski type equations for closed convex hypersurfaces were studied in \cite{DingLi2025}.
These developments show that the dual Orlicz setting is a natural nonhomogeneous counterpart of the dual \(L_p\) theory,
with the scaling behavior encoded by the radial Orlicz function.

In the capillary setting, one studies convex hypersurfaces in the half-space which meet the boundary hyperplane
at a fixed contact angle.
Recent work has developed the capillary analogue of the Christoffel--Minkowski theory,
including capillary Alexandrov--Fenchel inequalities, capillary support functions, curvature measures,
non-collapsing estimates and degree-theoretic existence methods;
see, for instance, \cite{HY2026,LiLiu2026,MWW2025,MWX2024,WWX2022,Zhu2026}.
Under the capillary Gauss parametrization, the geometric problem becomes a fully nonlinear equation on a
spherical cap with a Robin boundary condition.

Let \(0<\theta<\pi/2\) and define the normal spherical cap by
\[
 \Om=\C_\theta
 :=\{\xi=(\xi',\xi_{n+1})\in\Sph^n:\xi_{n+1}>\cos\theta\}.
\]
We denote the standard spherical metric by \(\sig\) and the outer unit conormal of \(\partial\Om\) by \(\mu\).
A function \(v\) on \(\overline\Om\) is called \emph{horizontally even} if
\[
 v(-\xi',\xi_{n+1})=v(\xi',\xi_{n+1}),
 \qquad \xi\in\overline\Om.
\]
For a positive support function \(h\), set
\[
 A[h]=\nabla^2h+h\sig,
 \qquad
 \rho[h]=\sqrt{h^2+|\nabla h|^2}.
\]
Here \(\sigma_k\) denotes the unnormalized \(k\)-th elementary symmetric polynomial in the eigenvalues:
\[
 \sigma_k(\lambda)
 =\sum_{1\le i_1<\cdots<i_k\le n}\lambda_{i_1}\cdots\lambda_{i_k},
 \qquad \sigma_k(\Id)=\binom nk.
\]
We consider the problem
\begin{equation}\label{eq:intro-main}
 \sigma_k(A[h])
 =
 f(\xi)h^{p-1}\frac{\rho[h]^{k+1}}{\varphi(\rho[h])}
 \quad\text{in }\Om,
\end{equation}
subject to
\begin{equation}\label{eq:intro-robin}
 \nabla_\mu h=\cot\theta\,h
 \quad\text{on }\partial\Om.
\end{equation}
When \(\varphi(r)=r^q\), equation \eqref{eq:intro-main} reduces to the capillary \(L_p\) dual
Christoffel--Minkowski equation with radial exponent \(q\).

The power-type equation therefore suggests replacing the homogeneous radial factor by a general Orlicz function
\(\varphi\). This leads naturally to \eqref{eq:intro-main}. We next explain its geometric interpretation through the
corresponding capillary dual Orlicz curvature measure.

The terminology ``dual Orlicz'' comes from the radial quantity in the capillary support parametrization.
If
\[
 X(\xi)=\nabla h(\xi)+h(\xi)\xi,
\]
then
\[
 |X(\xi)|
 =
 \sqrt{h^2+|\nabla h|^2}
 =
 \rho[h](\xi).
\]
Motivated by the dual Orlicz curvature measure introduced in \cite{ZXY2018}, we define in the smooth capillary setting
the \(L_p\)-weighted dual Orlicz \(k\)-th curvature measure by
\[
 d\widetilde C_{p,k,\varphi}^{\,\theta}(h,\xi)
 :=
 h^{1-p}
 \frac{\varphi(\rho[h])}{\rho[h]^{k+1}}
 \sigma_k(A[h])\,d\sigma(\xi).
\]
Consequently, prescribing
\[
 d\widetilde C_{p,k,\varphi}^{\,\theta}(h,\xi)
 =
 f(\xi)\,d\sigma(\xi)
\]
is exactly equivalent to \eqref{eq:intro-main}.
When \(\varphi(r)=r^q\), this reduces to the power-type capillary \(L_p\) dual Christoffel--Minkowski density
considered in \cite{LiLiu2026}; the corresponding closed problem was studied in \cite{CabezasHu2025}.

The purpose of the present paper is to identify natural hypotheses under which the boundary value problem
\eqref{eq:intro-main}--\eqref{eq:intro-robin} is solvable.
Two quantities associated with \(\varphi\) play distinct roles:
\[
 H(r)=\frac{\varphi(r)}{r^p},
 \qquad
 Q_\varphi(r)=\frac{r\varphi'(r)}{\varphi(r)}.
\]
The identity
\[
 \frac{rH'(r)}{H(r)}=Q_\varphi(r)-p
\]
shows that the strict increase of \(H\) is precisely the natural scale monotonicity.
The properness conditions
\[
 H(0+)=0,\qquad H(\infty)=\infty
\]
are used to rule out scale collapse and scale blow-up.
On the other hand, the upper bound
\[
 Q_\varphi\le k+1
\]
is the sign condition that enters the weighted gradient estimate.

The density condition is most naturally expressed through
\[
 g=f^{-1/(p+k-1)}.
\]
The inequalities
\[
 \nabla^2g+g\sig\ge0,
 \qquad
 \nabla_\mu g\le \cot\theta\,g
\]
are exactly the capillary spherical-convexity and boundary compatibility conditions appearing in the
constant-rank mechanism when \(k<n\).

Our main result is the following.

\begin{theorem}\label{thm:main}
Let
\[
 2\le k\le n,\qquad 1<p<k+1,\qquad 0<\theta<\frac{\pi}{2}.
\]
Let \(f\in C^\infty(\overline\Om)\) be positive and horizontally even.
Set
\[
 g=f^{-1/(p+k-1)}.
\]
Assume
\begin{equation}\label{eq:gconv}
 \nabla^2g+g\sig\ge0
 \quad\text{in }\Om,
\end{equation}
and
\begin{equation}\label{eq:gbdry}
 \nabla_\mu g\le\cot\theta\,g
 \quad\text{on }\partial\Om.
\end{equation}
Let \(\varphi\in C^\infty((0,\infty))\) be positive.
Assume that
\begin{equation}\label{eq:Hcondition}
 H(r):=\frac{\varphi(r)}{r^p}
\end{equation}
is strictly increasing and satisfies
\begin{equation}\label{eq:Hproper}
 \lim_{r\downarrow0}H(r)=0,
 \qquad
 \lim_{r\uparrow\infty}H(r)=\infty,
\end{equation}
and that
\begin{equation}\label{eq:Qupper}
 Q_\varphi(r):=\frac{r\varphi'(r)}{\varphi(r)}\le k+1,
 \qquad r>0.
\end{equation}
Then the boundary value problem \eqref{eq:intro-main}--\eqref{eq:intro-robin} admits a smooth, positive,
horizontally even solution \(h\) such that
\[
 \nabla^2h+h\sig>0
 \quad\text{on }\overline\Om.
\]
\end{theorem}

The proof is based on the following ingredients.
First, we construct an Orlicz homotopy that preserves both the scale cone and the capillary spherical-convexity cone.
Second, we use the weighted test function of Hu and Ivaki \cite{HuIvaki2026}
and extend the power-type dual calculation of Li and Liu \cite{LiLiu2026}.
For the Orlicz weight, the coefficient \(k+1-q\) from differentiating the radial factor becomes
\(k+1-Q_\varphi(\rho)\), which is nonnegative.
Combining the resulting gradient estimate with the capillary John non-collapsing theorem for even
hypersurfaces gives uniform control of \(\max h/\min h\).
The properness of \(H\) then yields the two-sided \(C^0\) estimate.
Third, once the \(C^0\) and \(C^1\) estimates have been established, the Orlicz right-hand side has a uniform
\(C^2\) bound as a function of position, the solution value, and the gradient, treated as independent variables.
The global Hessian estimate can therefore be applied without circular dependence on \(D^2h\).
Fourth, for \(k<n\) we prove strict convexity by extending the capillary full-rank argument of
Hu, Ivaki and Scheuer \cite{HuIvakiScheuer2025} to the Orlicz equation.
In the interior calculation, the term involving \(Q_\varphi'\) is quadratic in the eigenvalue
associated with the differentiation direction, so it vanishes along a null eigenvector.
Finally, global oblique regularity and elliptic degree theory complete the continuation argument.

The paper is organized as follows.
In Sect.~2, we collect the capillary support-function identities, introduce the capillary dual Orlicz curvature measure,
and construct the Orlicz homotopy.
Sect.~3 is devoted to the a priori estimates, strict convexity, and higher regularity.
In Sect.~4, we complete the degree argument, prove \Cref{thm:main}, and present examples and remarks.

\section{Background, Notation and the Orlicz Homotopy}

\subsection{Capillary support and radial functions}

We work on the unit sphere \((\Sph^n,\sig)\).
The covariant derivative and Laplacian with respect to \(\sig\) are denoted by \(\nabla\) and \(\Delta\).
We use the normal spherical cap \(\Om=\C_\theta\) and its outer unit conormal \(\mu\), as defined in Sect.~1.
Here \(\xi\) denotes the unit normal. The translated capillary spherical cap used in \cite{LiLiu2026}
is identified with \(\overline\Om\) by the isometry
\[
 \xi\longmapsto\xi-\cos\theta\,e_{n+1},
\]
where \(e_{n+1}\) is the last coordinate unit vector in \(\mathbb R^{n+1}\).
Support functions and spherical derivatives are pulled back under this identification.

We fix the horizontally even reference function
\[
 \ell(\xi)=1-\cos\theta\,\xi_{n+1}.
\]
It is smooth and positive on \(\overline\Om\) and satisfies
\begin{equation}\label{eq:ell}
 \nabla^2\ell+\ell\sig=\sig,
 \qquad
 \nabla_\mu\ell=\cot\theta\,\ell.
\end{equation}
For \(h\in C^2(\overline\Om)\), define
\[
 A[h]=\nabla^2h+h\sig.
\]
For the geometric terminology below, a capillary support function is required to be positive and to satisfy
\(\nabla_\mu h=\cot\theta\,h\) on \(\partial\Om\).
Such a function is called \emph{capillary convex} if \(A[h]\ge0\), and
\emph{strictly capillary convex} if \(A[h]>0\), with the tensor inequalities understood on \(\overline\Om\).

For the differential identities below, assume that \(h\) is smooth and positive.
The tensor \(A[h]\) is Codazzi:
\begin{equation}\label{eq:codazzi}
 A_{ij;k}=A_{ik;j}.
\end{equation}
Set
\[
 \rho=\sqrt{h^2+|\nabla h|^2}.
\]
Under the capillary support parametrization
\begin{equation}\label{eq:support-param}
 X(\xi)=\nabla h(\xi)+h(\xi)\xi,
\end{equation}
the boundary identity
\[
 \langle X,e_{n+1}\rangle
 =h\cos\theta-\sin\theta\,\nabla_\mu h
 \quad\text{on }\partial\Om
\]
shows that the Robin condition places the boundary of the reconstructed hypersurface in the hyperplane
\(\{x_{n+1}=0\}\). Moreover, one has
\begin{equation}\label{eq:radial-length}
 |X(\xi)|^2
 =
 h(\xi)^2+|\nabla h(\xi)|^2,
\end{equation}
and hence
\[
 \rho(\xi)=|X(\xi)|.
\]
Thus \(\rho\) is the radial quantity naturally associated with the capillary support function.

A direct computation gives
\begin{equation}\label{eq:Rfirst}
 \rho\rho_i=h_jA_{ji},
\end{equation}
and
\begin{equation}\label{eq:Rsecond}
 \rho_{ij}
 =
 \frac{(A^2)_{ij}-hA_{ij}+h_mA_{ij;m}}\rho
 -
 \frac{(h_mA_{mi})(h_sA_{sj})}{\rho^3}.
\end{equation}

\subsection{The capillary dual Orlicz curvature measure}

In the classical dual Orlicz theory, Zhu, Xing and Ye \cite{ZXY2018} introduced the dual Orlicz curvature measure
by weighting the radial Gauss image with a general positive function of the radial function.
In the present smooth capillary setting, the radial quantity is \(\rho=|X|\) by \eqref{eq:radial-length}.
Accordingly, for a smooth positive strictly capillary-convex support function \(h\) satisfying the Robin condition,
we adopt the following convention for the capillary dual Orlicz \(k\)-th curvature measure:
\begin{equation}\label{eq:dual-orlicz-measure}
 d\widetilde C_{k,\varphi}^{\,\theta}(h,\xi)
 :=
 \frac{\varphi(\rho)}{\rho^{k+1}}
 \sigma_k(A[h])\,d\sigma(\xi).
\end{equation}
Its \(L_p\)-weighted version is
\begin{equation}\label{eq:lp-dual-orlicz-measure}
 d\widetilde C_{p,k,\varphi}^{\,\theta}(h,\xi)
 :=
 h^{1-p}
 \frac{\varphi(\rho)}{\rho^{k+1}}
 \sigma_k(A[h])\,d\sigma(\xi).
\end{equation}
For comparison, when \(k=n\), the classical dual Orlicz curvature measure of
\cite{ZXY2018}, expressed in the same normal coordinates, has the local density
\[
 \frac{1}{n+1}\,h\,
 \frac{\varphi(\rho)}{\rho^{n+1}}
 \sigma_n(A[h])\,d\sigma(\xi).
\]
Thus, in this top-order case, the density in \eqref{eq:dual-orlicz-measure}
is \((n+1)h^{-1}\) times the classical density, while that in
\eqref{eq:lp-dual-orlicz-measure} is \((n+1)h^{-p}\) times the classical density.
This comparison concerns the local densities under the common normal parametrization.

Therefore the smooth capillary dual Orlicz--Christoffel--Minkowski problem is to prescribe the density of
\(\widetilde C_{p,k,\varphi}^{\,\theta}\).
More precisely,
\begin{equation}\label{eq:measure-prescription}
 d\widetilde C_{p,k,\varphi}^{\,\theta}(h,\xi)
 =
 f(\xi)\,d\sigma(\xi)
\end{equation}
is equivalent to
\[
 \sigma_k(A[h])
 =
 f(\xi)h^{p-1}\frac{\rho^{k+1}}{\varphi(\rho)}.
\]

When
\[
 \varphi(r)=r^q,
\]
formula \eqref{eq:lp-dual-orlicz-measure} becomes
\[
 d\widetilde C_{p,k,q}^{\,\theta}(h,\xi)
 =
 h^{1-p}\rho^{q-k-1}\sigma_k(A[h])\,d\sigma(\xi),
\]
and the prescribed-measure equation becomes
\[
 \sigma_k(A[h])
 =
 f(\xi)h^{p-1}\rho^{k+1-q},
\]
which is precisely the power-type capillary \(L_p\) dual Christoffel--Minkowski equation treated in
\cite{LiLiu2026}.
Thus \eqref{eq:lp-dual-orlicz-measure} gives the natural smooth Orlicz extension of the power density.

On the cone of positive definite symmetric matrices, we shall use
\[
 S(A)=\sigma_k(A),
 \qquad
 F(A)=\log \sigma_k(A),
\]
and the matrix derivatives
\[
 S^{ij}=\frac{\partial S}{\partial A_{ij}},
 \qquad
 F^{ij}=\frac{\partial F}{\partial A_{ij}}=\frac{S^{ij}}{S}.
\]
If \(A\) is diagonal with eigenvalues
\[
 \lambda_1\ge\cdots\ge\lambda_n,
\]
then
\[
 \sum_iF^{ii}\lambda_i=k.
\]

For later use, set
\[
 G(r)=\frac{r^{k+1}}{\varphi(r)},
 \qquad
 a(r)=k+1-Q_\varphi(r).
\]
Then
\begin{equation}\label{eq:Gderivs}
 (\log G)'(r)=\frac{a(r)}r,
 \qquad
 (\log G)''(r)
 =
 -\frac{a(r)+rQ_\varphi'(r)}{r^2}.
\end{equation}
Combining \eqref{eq:Rfirst}, \eqref{eq:Rsecond} and \eqref{eq:Gderivs} yields
\begin{align}
 \nabla_{ij}\log G(\rho)
 ={}&
 \frac{a(\rho)}{\rho^2}
 \left((A^2)_{ij}-hA_{ij}+h_mA_{ij;m}\right)
 \notag\\
 &-
 \frac{2a(\rho)+\rho Q_\varphi'(\rho)}{\rho^4}
 (h_mA_{mi})(h_sA_{sj}).
 \label{eq:logGsecond}
\end{align}

\subsection{The Orlicz homotopy}

Let
\[
 d=p+k-1.
\]
At the reference endpoint, we choose
\[
 \varphi_0(r)=r^{k+1},
 \qquad
 f_0=\binom nk\ell^{1-p}.
\]
Then \(h_0=\ell\) solves the reference equation.

Set
\[
 g_0=f_0^{-1/d},
 \qquad
 g_1=f^{-1/d}.
\]
For \(t\in[0,1]\), define
\begin{equation}\label{eq:gt}
 g_t=(1-t)g_0+tg_1,
 \qquad
 f_t=g_t^{-d},
\end{equation}
and
\begin{equation}\label{eq:phit}
 \varphi_t(r)
 =
 r^{(k+1)(1-t)}\varphi(r)^t.
\end{equation}
We consider
\begin{equation}\label{eq:homotopy}
 \sigma_k(A[h])
 =
 f_t(\xi)h^{p-1}\frac{\rho^{k+1}}{\varphi_t(\rho)},
 \qquad
 \nabla_\mu h=\cot\theta\,h.
\end{equation}

The logarithmic slope satisfies
\begin{equation}\label{eq:Qt}
 Q_{\varphi_t}(r)
 =
 (1-t)(k+1)+tQ_\varphi(r).
\end{equation}
Since \(H\) is smooth and strictly increasing, \(H'(r)\ge0\). Thus
\[
 Q_\varphi(r)-p=\frac{rH'(r)}{H(r)}\ge0.
\]
Combining this with \eqref{eq:Qt} gives
\begin{equation}\label{eq:Qrange}
 p+(1-t)(k+1-p)\le Q_{\varphi_t}(r)\le k+1,
 \qquad r>0,\quad t\in[0,1].
\end{equation}
In particular, \(Q_{\varphi_t}(r)>p\) when \(t<1\), while equality with \(p\) is allowed at \(t=1\).
The coefficient used in the gradient estimate therefore satisfies
\[
 0\le k+1-Q_{\varphi_t}(r)
 =t\bigl(k+1-Q_\varphi(r)\bigr)
 \le t(k+1-p).
\]

Define
\[
 H_t(r)=\frac{\varphi_t(r)}{r^p}.
\]
Then
\begin{equation}\label{eq:Ht}
 H_t(r)
 =
 \bigl(r^{k+1-p}\bigr)^{1-t}H(r)^t.
\end{equation}

\begin{lemma}\label[lemma]{lem:proper}
The family \(\{H_t\}_{t\in[0,1]}\) is uniformly proper:
\[
 \sup_{t\in[0,1]}H_t(r)\to0
 \quad\text{as }r\downarrow0,
\]
and
\[
 \inf_{t\in[0,1]}H_t(r)\to\infty
 \quad\text{as }r\uparrow\infty.
\]
\end{lemma}

\begin{proof}
By \eqref{eq:Ht},
\[
 \min\{r^{k+1-p},H(r)\}\le H_t(r)
 \le\max\{r^{k+1-p},H(r)\},
 \qquad t\in[0,1].
\]
Both endpoint functions tend to zero as \(r\downarrow0\) and to infinity as \(r\uparrow\infty\),
which proves the stated uniform limits.
\end{proof}

To check the reference endpoint, write
\[
 g_0=b_0\ell^\alpha,
 \qquad b_0=\bigl(\tbinom nk\bigr)^{-1/d},
 \qquad \alpha=\frac{p-1}{p+k-1}\in(0,1).
\]
Since \(A[\ell]=\sig\) and \(|\nabla\ell|^2\le\cos^2\theta\,\ell\), direct differentiation gives
\[
 A[g_0]\ge
 b_0\ell^{\alpha-1}
 \bigl[\alpha\bigl(1-(1-\alpha)\cos^2\theta\bigr)+(1-\alpha)\ell\bigr]\sig>0.
\]
Moreover, \(\nabla_\mu g_0=\alpha\cot\theta\,g_0\le\cot\theta\,g_0\) on \(\partial\Om\).
The capillary spherical-convexity cone is therefore preserved:
\[
 \nabla^2g_t+g_t\sig
 =
 (1-t)(\nabla^2g_0+g_0\sig)
 +
 t(\nabla^2g_1+g_1\sig)\ge0,
\]
and
\[
 \nabla_\mu g_t-\cot\theta\,g_t\le0.
\]

\section{A Priori Estimates and Strict Convexity}

\subsection{The \(C^0\) and \(C^1\) Estimates}

\subsubsection*{Weighted gradient estimate}

Set
\[
 u=\frac h\ell.
\]
By the common Robin boundary condition,
\[
 \nabla_\mu u
 =
 \frac{\ell\nabla_\mu h-h\nabla_\mu\ell}{\ell^2}
 =0
 \quad\text{on }\partial\Om.
\]
Fix
\[
 0<\gamma<
 \min\left\{\frac{2(p-1)}n,2\right\}.
\]
Following Hu and Ivaki \cite[Lemma~3.1]{HuIvaki2026}, we use the capillary weighted test function
\[
 \ell^{2-\gamma}\frac{|\nabla u|^2}{u^\gamma}.
\]

\begin{proposition}\label[proposition]{prop:gradient}
There exists a constant \(C>0\), independent of \(t\), such that every smooth positive strictly capillary-convex solution of
\eqref{eq:homotopy} satisfies
\begin{equation}\label{eq:gradient}
 \sup_{\overline\Om}
 \frac{|\nabla h|^2}{h^\gamma}
 \le
 C\bigl(\max_{\overline\Om}h\bigr)^{2-\gamma}.
\end{equation}
\end{proposition}

\begin{proof}
To adapt the power-type dual calculation of Li and Liu \cite[Lemma~3.1]{LiLiu2026},
differentiate the logarithm of the radial factor in \eqref{eq:homotopy}. The chain rule gives
\[
 \nabla\log\left(\frac{\rho^{k+1}}{\varphi_t(\rho)}\right)
 =
 a_t(\rho)\frac{\nabla\rho}{\rho},
 \qquad
 a_t(r)=k+1-Q_{\varphi_t}(r).
\]
Thus \(a_t(\rho)\) takes the place of \(k+1-q\) in the differentiated equation.
If the desired bound fails, the ratio \(|\nabla u|/u\) is large at a maximum of the test function.
The estimate following (3.16) in \cite{LiLiu2026} shows that the factor multiplying \(a_t(\rho)\)
is nonnegative there.
The radial contribution is therefore nonnegative by \eqref{eq:Qrange}, including when \(a_t(\rho)=0\).
The estimate of the remaining terms uses the positivity of the radial factor and a uniform bound on
\(\nabla\log f_t=-d\,\nabla g_t/g_t\), which follows from \eqref{eq:gt}.
For \(|\nabla u|/u\) sufficiently large, the sum of those terms is strictly positive, contradicting
the maximum-principle inequality. This proves \eqref{eq:gradient} with a constant independent of \(t\).
\end{proof}

\subsubsection*{Non-collapsing and the \(C^0\) estimate}

The non-collapsing estimate of Hu and Yang \cite[Lemma~4.1]{HY2026} applies to horizontally even,
smooth, strictly convex capillary hypersurfaces satisfying \eqref{eq:gradient}. Thus
\begin{equation}\label{eq:ratio}
 \frac{\max_{\overline\Om}h}{\min_{\overline\Om}h}\le C.
\end{equation}

Let
\[
 M=\max_{\overline\Om}\frac h\ell,
 \qquad
 m=\min_{\overline\Om}\frac h\ell.
\]
Let \(x_M\in\overline\Om\) be a maximum point of \(h/\ell\).
At \(x_M\),
\[
 \nabla h=M\nabla\ell,
 \qquad
 \nabla^2h\le M\nabla^2\ell.
\]
Therefore, using \eqref{eq:ell},
\[
 A[h]
 =
 \nabla^2h+h\sig
 \le
 M(\nabla^2\ell+\ell\sig)
 =
 M\sig.
\]
Consequently,
\begin{equation}\label{eq:sigmakmax}
 \sigma_k(A[h](x_M))
 \le
 \binom nk M^k.
\end{equation}
Moreover,
\[
 h(x_M)=M\ell(x_M),
 \qquad
 \nabla h(x_M)=M\nabla\ell(x_M),
\]
and hence
\[
 \rho[h](x_M)
 =
 M \rho[\ell](x_M),
\]
where
\[
 \rho[\ell]=\sqrt{\ell^2+|\nabla\ell|^2}.
\]
Substituting these relations into \eqref{eq:homotopy} and using \eqref{eq:sigmakmax}, we obtain
\[
 \binom nk M^k
 \ge
 f_t(x_M)(M\ell)^{p-1}
 \frac{(M\rho[\ell])^{k+1}}{\varphi_t(M\rho[\ell])}.
\]
Since
\[
 \varphi_t(r)=r^pH_t(r),
\]
this becomes
\[
 \binom nk
 \ge
 f_t(x_M)\ell(x_M)^{p-1}
 \rho[\ell](x_M)^{k+1-p}
 \frac1{H_t(M\rho[\ell](x_M))}.
\]
Thus
\begin{equation}\label{eq:scale-upper}
 H_t(M\rho[\ell](x_M))
 \ge
 \frac{
 f_t(x_M)\ell(x_M)^{p-1}\rho[\ell](x_M)^{k+1-p}
 }{\binom nk}.
\end{equation}
The fixed reference function \(\ell\), the function \(\rho[\ell]\), and the homotopy densities \(f_t\) are uniformly
positive and bounded on \(\overline\Om\).
Therefore there exists \(c_0>0\), independent of \(t\), such that
\begin{equation}\label{eq:scale-upper2}
 H_t(M\rho[\ell](x_M))\ge c_0.
\end{equation}

Similarly, let \(x_m\in\overline\Om\) be a minimum point of \(h/\ell\).
Then
\[
 \nabla h=m\nabla\ell,
 \qquad
 \nabla^2h\ge m\nabla^2\ell,
\]
and hence
\[
 A[h](x_m)\ge m\sig.
\]
Therefore
\[
 \sigma_k(A[h](x_m))
 \ge
 \binom nk m^k.
\]
Since
\[
 h(x_m)=m\ell(x_m),
 \qquad
 \rho[h](x_m)=m\rho[\ell](x_m),
\]
the equation gives
\[
 \binom nk m^k
 \le
 f_t(x_m)(m\ell)^{p-1}
 \frac{(m\rho[\ell])^{k+1}}{\varphi_t(m\rho[\ell])}.
\]
Equivalently,
\begin{equation}\label{eq:scale-lower}
 H_t(m\rho[\ell](x_m))
 \le
 \frac{
 f_t(x_m)\ell(x_m)^{p-1}\rho[\ell](x_m)^{k+1-p}
 }{\binom nk}
 \le C_0,
\end{equation}
for some constant \(C_0>0\) independent of \(t\).

By the uniform properness in \cref{lem:proper}, \eqref{eq:scale-upper2} excludes
\[
 M\to0,
\]
and therefore
\begin{equation}\label{eq:Mlower}
 M\ge c>0.
\end{equation}
Likewise, \eqref{eq:scale-lower} excludes
\[
 m\to\infty,
\]
and hence
\begin{equation}\label{eq:mupper}
 m\le C.
\end{equation}

Because \(\ell\) is positive and smooth on \(\overline\Om\), there exist constants
\(0<c_\ell\le C_\ell\) such that
\[
 c_\ell\le\ell\le C_\ell.
\]
Therefore
\[
 \max_{\overline\Om}h
 \le
 C_\ell M,
 \qquad
 \min_{\overline\Om}h
 \ge
 c_\ell m,
\]
and also
\[
 \max_{\overline\Om}h
 \ge
 c_\ell M,
 \qquad
 \min_{\overline\Om}h
 \le
 C_\ell m.
\]
Combining these comparisons with \eqref{eq:ratio}, \eqref{eq:Mlower} and \eqref{eq:mupper}, we obtain
\[
 M\le C,
 \qquad
 m\ge c.
\]
Consequently,
\begin{equation}\label{eq:C0}
 0<c\le h\le C
 \quad\text{on }\overline\Om,
\end{equation}
with constants independent of \(t\).

\subsubsection*{The \(C^1\) estimate}

Once \eqref{eq:C0} is known, the weighted gradient estimate \eqref{eq:gradient} gives
\[
 |\nabla h|^2
 \le
 C h^\gamma
 \bigl(\max_{\overline\Om}h\bigr)^{2-\gamma}
 \le C.
\]
Thus
\begin{equation}\label{eq:C1}
 |\nabla h|\le C
 \quad\text{on }\overline\Om.
\end{equation}

We have therefore proved the following.

\begin{proposition}\label[proposition]{prop:c01}
There exists a constant \(C\ge1\), independent of \(t\), such that every positive, horizontally even, smooth, strictly capillary-convex solution
of \eqref{eq:homotopy} satisfies
\[
 C^{-1}\le h\le C,
 \qquad
 |\nabla h|\le C.
\]
\end{proposition}

\subsection{The Global Hessian Estimate}

Define the low-order function
\begin{equation}\label{eq:fhat}
 \widehat f_t(\xi,z,q)
 =
 f_t(\xi)z^{p-1}
 \frac{(z^2+|q|^2)^{(k+1)/2}}
 {\varphi_t\bigl((z^2+|q|^2)^{1/2}\bigr)}.
\end{equation}
The values \((z,q)=(h,\nabla h)\) generated by the homotopy solutions covered by \cref{prop:c01}
lie in a fixed compact subset of \((0,\infty)\times\mathbb R^n\).
The derivatives in the following lemma are taken with respect to the independent variables \((\xi,z,q)\).

\begin{lemma}\label[lemma]{lem:fhat}
There exists \(C>0\), independent of \(t\), such that
\[
 C^{-1}\le\widehat f_t\le C,
 \qquad
 \|\widehat f_t\|_{C^2}\le C
\]
on a fixed compact low-order range containing these values.
\end{lemma}

\begin{proof}
Let
\[
 G_t(r)=\frac{r^{k+1}}{\varphi_t(r)},
 \qquad
 a_t(r)=k+1-Q_{\varphi_t}(r).
\]
Then
\[
 \frac{G_t'}{G_t}=\frac{a_t(r)}r,
\]
and
\[
 \frac{G_t''}{G_t}
 =
 \frac{a_t(r)^2-a_t(r)-rQ_{\varphi_t}'(r)}{r^2}.
\]
On every fixed compact subinterval of \((0,\infty)\), smoothness and positivity of \(\varphi\) imply uniform bounds for
\(\varphi_t\), \(\varphi_t^{-1}\), \(Q_{\varphi_t}\) and \(rQ_{\varphi_t}'\).
The claim follows by the chain rule.
\end{proof}

The global \(C^2\) estimate in \cite[Theorem~3.5]{LiLiu2026} is formulated for
\[
 \sigma_k(\nabla^2h+h\sig)=\widehat f(\xi,h,\nabla h)
\]
with a positive right-hand side whose \(C^2\)-norm in the independent lower-order variables is controlled.
Thus \cref{lem:fhat} provides the corresponding Orlicz verification.

For later reference, define
\[
 E_\ell=\frac12\ell^2\left|\nabla\left(\frac h\ell\right)\right|^2.
\]
In an adapted orthonormal frame, with \(\alpha\) denoting a direction tangent to \(\partial\Om\),
the condition \(u_\mu=0\) for \(u=h/\ell\) gives
\[
 u_{\alpha\mu}=-\cot\theta\,u_\alpha,
\]
and hence
\begin{equation}\label{eq:Ebdry}
 (E_\ell)_\mu=0
 \quad\text{on }\partial\Om.
\end{equation}
This is the key boundary identity used in the power-sum Hessian auxiliary function.

\begin{proposition}\label[proposition]{prop:c2}
There exists \(C>0\), independent of \(t\), such that every positive, horizontally even, smooth,
strictly capillary-convex solution of \eqref{eq:homotopy} satisfies
\[
 \|h\|_{C^2(\overline\Om)}\le C.
\]
Equivalently,
\[
 0<A[h]\le C\sig.
\]
\end{proposition}

\begin{proof}
By \cref{prop:c01}, there are constants \(c_0,C_0,C_1>0\), independent of \(t\), such that
\[
 c_0\le h\le C_0,
 \qquad
 |\nabla h|\le C_1.
\]
Hence
\[
 c_0\le \rho=\sqrt{h^2+|\nabla h|^2}\le C
\]
for a uniform constant \(C\). Thus the arguments of \(\varphi_t\) stay in a fixed compact interval contained in
\((0,\infty)\). By \cref{lem:fhat}, the function \(\widehat f_t(\xi,z,q)\), with \((\xi,z,q)\)
treated as independent variables, is uniformly positive and has a uniformly bounded \(C^2\)-norm
on a fixed compact range containing the values \((\xi,h,\nabla h)\) generated by these solutions.

We may therefore write the equation as
\[
 \sigma_k(\nabla^2h+h\sig)=\widehat f_t(\xi,h,\nabla h)
\]
and apply the global Hessian estimate in \cite[Theorem~3.5]{LiLiu2026}.
We record the two points of that argument relevant to the present equation.

First, the boundary test function contains
\[
 E_\ell
 =
 \frac12\ell^2\left|\nabla\left(\frac h\ell\right)\right|^2.
\]
Because \(u=h/\ell\) satisfies \(u_\mu=0\), the capillary boundary identities give
\[
 (E_\ell)_\mu=0
\]
on \(\partial\Om\), as in \eqref{eq:Ebdry}. Thus the large coefficient multiplying \(E_\ell\) does not contribute to the
boundary normal derivative estimate. The remaining boundary terms are handled as in the proof of
\cite[Theorem~3.5]{LiLiu2026}, excluding a boundary maximum of the auxiliary function uniformly in \(t\).

Second, at an interior maximum of the Hessian auxiliary function, diagonalize \(A[h]\) and write
\[
 \lambda_1\ge\cdots\ge\lambda_n>0.
\]
We use the parameter choices in the proof of
\cite[Theorem~3.5]{LiLiu2026}.  In that proof, the coefficient of
\(E_\ell\) is denoted by \(A\); we denote it by \(\Lambda\) here to
distinguish it from \(A[h]\).  For \(k<n\), let
\(\delta=\delta_k\in(0,1)\) be the final spectral threshold constructed in
\cite[Lemma~3.8]{LiLiu2026}.  These thresholds and the auxiliary parameters
\(W,K,L\) in that proof are fixed before \(\Lambda\), using the uniform
data bounds established above; the boundary choice is independent of
\(\Lambda\) because of \eqref{eq:Ebdry}.

In the spectral-gap branch, the third-derivative estimates and the
maximum-principle inequality give, as in
\cite[(3.42)--(3.46)]{LiLiu2026},
\[
 0\ge
 \Lambda\sigma_k^{11}\lambda_1^2-C_1\lambda_1-C_0
 -\frac{C_{-1}}{\lambda_1}.
\]
Here \(C_1,C_{-1}\) are independent of \(\Lambda\), while \(C_0\) may
depend on \(\Lambda\) and the subsequently chosen auxiliary parameters.
This dependence is the one stated after (3.38) in the cited proof.
Since \(\lambda_1\) is the largest eigenvalue,
\[
 \lambda_1\sigma_k^{11}
 \ge\frac{k}{n}\sigma_k(A[h]),
 \qquad
 \sigma_k^{11}\lambda_1^2\ge c_*\lambda_1,
 \qquad
 c_*=\frac{k}{n}s_0>0,
\]
where \(s_0\) is the uniform positive lower bound for
\(\sigma_k(A[h])\).  Choose \(\Lambda\) large enough to satisfy the
auxiliary-parameter hierarchy in that proof and
\(\Lambda c_*-C_1\ge1\), and then choose the remaining parameters as
there.  The inequality becomes
\[
 0\ge\lambda_1-C_0-\frac{C_{-1}}{\lambda_1},
\]
which bounds \(\lambda_1\) uniformly.  In the complementary branch,
\(\lambda_k\ge\delta\lambda_1\), and hence
\[
 \sigma_k(A[h])\ge\lambda_1\cdots\lambda_k
 \ge\delta^{k-1}\lambda_1^k.
\]
The uniform upper bound for \(\sigma_k(A[h])\) again bounds
\(\lambda_1\).

For \(k=n\), let \(\delta=\delta_n\) be the corresponding final threshold
from the cited spectral argument.  The same gap-branch inequality closes using
\[
 \sigma_n^{11}\lambda_1^2=\det(A[h])\lambda_1
 \ge s_0\lambda_1.
\]
In the complementary branch, the final spectral threshold from the same
argument gives \(\lambda_n\ge\delta\lambda_1\), so
\[
 \det A[h]\ge\delta^{n-1}\lambda_1^n,
\]
and the uniform upper bound for \(\det A[h]\) yields
\(\lambda_1\le C\).  All parameter choices and resulting bounds are
uniform in \(t\), since the data controlling the cited estimate are
uniform in \(t\).

Thus the largest eigenvalue of \(A[h]\) is uniformly bounded on \(\overline\Om\), and consequently
\[
 \|h\|_{C^2(\overline\Om)}\le C.
\]
\end{proof}

\subsection{Constant Rank and Strict Convexity}

We adapt the capillary full-rank argument of Hu, Ivaki and Scheuer \cite[Theorem~3.1]{HuIvakiScheuer2025}
to \eqref{eq:intro-main} in the case \(k<n\).
We verify the interior viscosity structural inequality for the Orlicz right-hand side.
At the boundary we follow their analysis of tangential and normal null directions.
Throughout this subsection, \(A=A[h]=\nabla^2h+h\sig\).
Let \(h\in C^\infty(\overline\Om)\) be a positive solution of
\eqref{eq:intro-main}--\eqref{eq:intro-robin} with \(A[h]\ge0\).
We assume \eqref{eq:gconv}, and additionally \eqref{eq:gbdry} when addressing the boundary.
The spherical curvature identity implies that \(A\) is a Codazzi tensor.

Set
\[
 \Psi(A)=\sigma_k(A)^{1/k},
\]
and
\[
 \Hh
 =
 \left(
 f h^{p-1}\frac{\rho^{k+1}}{\varphi(\rho)}
 \right)^{1/k}.
\]
Then
\[
 \Psi(A)=\Hh.
\]
Here the G\aa rding cone is
\[
 \Gamma_k
 =\{\lambda\in\mathbb R^n:\sigma_j(\lambda)>0\text{ for }1\le j\le k\},
\]
and the same notation denotes its associated spectral matrix domain.
The function \(\Psi\) is smooth, positive, one-homogeneous and strictly elliptic on \(\Gamma_k\).
It is inverse-concave in the convention of \cite[Definitions~1.3--1.4]{BIS2020}:
the inverse-matrix test is made on the positive definite cone \(\Gamma_+\subset\Gamma_k\).
In particular, at a positive definite matrix \(\beta\), for every symmetric tensor \(\eta\),
\[
 \Psi^{ij,rs}(\beta)\eta_{ij}\eta_{rs}
 +2\Psi^{ir}(\beta)(\beta^{-1})^{js}\eta_{ij}\eta_{rs}
 \ge \frac{2}{\Psi(\beta)}\bigl(\Psi^{ij}(\beta)\eta_{ij}\bigr)^2;
\]
see \cite[(2.1) and Section~2]{BIS2020}.
The derivatives of \(\Psi\) are smooth on \(\Gamma_k\), including at its nonnegative matrices
with \(\sigma_k>0\); the displayed inverse-matrix test is needed only at positive definite matrices.

Let
\[
 \tau_0=\frac{p-1}{k},
 \qquad
 d=p+k-1=k(1+\tau_0),
 \qquad
 g=f^{-1/d},
\]
and set
\[
 W=\Hh^{-1}.
\]
Writing
\[
 G(\rho)=\frac{\rho^{k+1}}{\varphi(\rho)},
\]
we have
\begin{equation}\label{eq:W}
 W=g^{1+\tau_0}h^{-\tau_0}G(\rho)^{-1/k}.
\end{equation}

Fix a point and choose an \(A\)-diagonal orthonormal frame.
Let \(e_m\) be an eigenvector with eigenvalue \(\lambda\).
By \eqref{eq:Rfirst}--\eqref{eq:Rsecond},
\[
 \rho_m=\frac{h_m\lambda}{\rho},
\]
and
\[
 \rho_{mm}
 =
 \frac{\lambda^2-h\lambda+h_jA_{mm;j}}\rho
 -
 \frac{h_m^2\lambda^2}{\rho^3}.
\]

\begin{lemma}\label[lemma]{lem:Winv}
Let \(a(\rho)=k+1-Q_\varphi(\rho)\).
At the above point,
\begin{align}
 \frac{W_{mm}+W}{W}
 ={}&
 (1+\tau_0)\frac{g_{mm}+g}{g}
 +
 \tau_0(1+\tau_0)
 \left(
 \frac{g_m}{g}-\frac{h_m}{h}
 \right)^2
 \notag\\
 &-
 \frac{a(\rho)}{k\rho^2}h_jA_{mm;j}
 +
 \mathcal E_m,
 \label{eq:Winv}
\end{align}
where every term in \(\mathcal E_m\) contains at least one factor \(\lambda\), and the unique term containing
\(Q_\varphi'(\rho)\) is proportional to \(\lambda^2Q_\varphi'(\rho)\).
Consequently, on every compact subset of \(\Om\), with \(C\) depending on the fixed solution,
\begin{equation}\label{eq:Winvbound}
 W_{mm}+W
 \ge
 -C\left(\lambda+|\nabla A_{mm}|\right).
\end{equation}
\end{lemma}

\begin{proof}
Differentiate \eqref{eq:W} twice.
Using
\[
 (\log G)'(\rho)=\frac{a(\rho)}\rho,
 \qquad
 (\log G)''(\rho)
 =
 -\frac{a(\rho)+\rho Q_\varphi'(\rho)}{\rho^2},
\]
one obtains
\begin{align*}
 \mathcal E_m
 ={}&
 -\frac{\tau_0\lambda}{h}
 -\frac{a(\rho)}{k\rho^2}(\lambda^2-h\lambda)
 +
 \frac{2a(\rho)+\rho Q_\varphi'(\rho)}{k\rho^4}h_m^2\lambda^2
\\
 &-
 \frac{2a(\rho)h_m\lambda}{k\rho^2}
 \left(
 (1+\tau_0)\frac{g_m}{g}
 -
 \tau_0\frac{h_m}{h}
 \right)
 +
 \frac{a(\rho)^2h_m^2\lambda^2}{k^2\rho^4}.
\end{align*}
The first two terms on the right side of \eqref{eq:Winv} are nonnegative because
\[
 g_{mm}+g\ge0,
 \qquad
 \tau_0\ge0.
\]
The remaining terms give \eqref{eq:Winvbound}.
\end{proof}

For the top-order case \(k=n\), full rank is immediate. Indeed,
\[
 \det A[h]
 =
 f h^{p-1}\frac{\rho^{n+1}}{\varphi(\rho)}
 >0.
\]
Since admissibility gives \(A[h]\ge0\), every eigenvalue of \(A[h]\) is positive, and hence
\[
 A[h]>0
 \quad\text{on }\overline\Om.
\]
Thus, in the remainder of the constant-rank argument, we assume \(k<n\).

For the interior step, we apply the viscosity constant-rank theorem of Bryan, Ivaki and Scheuer
\cite[Theorem~1.6]{BIS2020}. To verify its structural hypothesis, define
\[
 \mathscr F(A,x)=\Psi(A)-\Hh(x),
 \qquad
 \Phi^{ij,rs}(A,x)=\frac{2}{\Psi(A)}\Psi^{ij}(A)\Psi^{rs}(A).
\]
Here \(\Hh(x)\) is the smooth scalar field determined by the fixed solution \(h\).
Its derivatives below are total derivatives along that solution.
When taking the explicit position derivatives in \(\mathscr F\), the matrix argument is held fixed;
thus \(D_x\mathscr F^{ij}=0\) and \(D^2_{xx}\mathscr F=-\nabla^2\Hh\).
The preceding inverse-concavity inequality is precisely the
\(\Phi\)-inverse concavity required in \cite[Definition~1.4]{BIS2020}.

At a point choose an \(A\)-diagonal orthonormal frame and a unit eigenvector \(e_m\),
with eigenvalue \(\lambda=A_{mm}\).
The curvature-adjusted modulus \(\omega_{\mathscr F}\) of that definition, specialized to the unit sphere,
is given for a symmetric tensor \(\eta\) by
\[
 \omega_{\mathscr F}(A)(\eta,e_m)
 =\frac{2}{\Psi(A)}\bigl(\Psi^{ij}\eta_{ij}\bigr)^2
  -\Hh_{mm}+\sum_{i\ne m}\lambda_i\Psi^{ii}.
\]
Differentiating \(\Psi(A)=\Hh\) gives
\[
 \Psi^{ij}A_{ij;m}=\Hh_m,
\]
and one-homogeneity gives
\[
 \sum_i\Psi^{ii}\lambda_i=\Psi(A)=\Hh,
 \qquad
 \sum_{i\ne m}\lambda_i\Psi^{ii}=\Hh-\lambda\Psi^{mm}.
\]
Consequently the exact modulus identity is
\[
 \omega_{\mathscr F}(A)(\nabla_m A,e_m)
 =\frac{2\Hh_m^2}{\Hh}-\Hh_{mm}+\Hh-\lambda\Psi^{mm}
 =\Hh^2(W_{mm}+W)-\lambda\Psi^{mm}.
\]
By \eqref{eq:Winvbound}, on each compact subset of \(\Om\),
\[
 \omega_{\mathscr F}(A)(\nabla_m A,e_m)
 \ge -C\bigl(\lambda+|\nabla A_{mm}|\bigr).
\]
This is the eigenvector structural inequality in \cite[Theorem~1.6]{BIS2020}.
Since \(A\ge0\) and \(\sigma_k(A)=\Hh^k>0\), its eigenvalues lie in \(\Gamma_k\).
All hypotheses of that theorem are therefore satisfied, and \(A\) has constant rank
in the connected interior of \(\Om\).

To identify this rank, note that
\[
 \nabla_\mu h=\cot\theta\,h>0
 \quad\text{on }\partial\Om.
\]
Thus a minimum of \(h\) cannot occur on the boundary, and \(h\) attains its minimum at an interior point \(x_0\).
At \(x_0\),
\[
 \nabla h(x_0)=0,
 \qquad
 \nabla^2h(x_0)\ge0.
\]
Therefore
\[
 A[h](x_0)
 =
 \nabla^2h(x_0)+h(x_0)\sig
 \ge
 h(x_0)\sig>0.
\]
Hence \(\operatorname{rank}A(x_0)=n\). By constancy of the rank,
\[
 A[h]>0
 \quad\text{in }\Om.
\]

\begin{proposition}\label[proposition]{prop:interior-full}
Let \(2\le k\le n\), \(1<p<k+1\), and \(0<\theta<\pi/2\).
Let \(f\in C^\infty(\overline\Om)\) and \(\varphi\in C^\infty((0,\infty))\) be positive.
Suppose that \(h\in C^\infty(\overline\Om)\) is a positive solution of
\eqref{eq:intro-main}--\eqref{eq:intro-robin}, that \(A[h]\ge0\), and that
\(g=f^{-1/(p+k-1)}\) satisfies \(\nabla^2g+g\sig\ge0\). Then
\[
 A[h]>0
 \quad\text{in }\Om.
\]
\end{proposition}

It remains to rule out rank loss at the boundary.
Fix \(x_0\in\partial\Om\) and choose an orthonormal frame with
\[
 e_n=\mu,
\]
while \(e_\alpha\), \(1\le\alpha\le n-1\), are tangent to \(\partial\Om\),
and choose them to diagonalize the tangential block of \(A(x_0)\).
The Robin condition and the capillary boundary geometry give
\begin{equation}\label{eq:block}
 A_{\alpha n}=0.
\end{equation}
Hu, Ivaki and Scheuer's boundary derivative identity \cite[Eq.~(3.3)]{HuIvakiScheuer2025} gives
\begin{equation}\label{eq:bdryA}
 A_{\alpha\beta;n}
 =
 \cot\theta
 \left(
 A_{nn}\delta_{\alpha\beta}-A_{\alpha\beta}
 \right).
\end{equation}

By \eqref{eq:block}, if \(A(x_0)\) has a nontrivial kernel, it contains either a tangential
unit null vector or the normal vector \(\mu\). It is enough to consider these two alternatives.
Suppose first that a tangential direction \(e_\alpha\) is null. Then \(A_{\alpha\alpha}=0\), and
\eqref{eq:bdryA} gives
\[
 A_{\alpha\alpha;n}
 =
 \cot\theta\,A_{nn}\ge0.
\]
Thus the tensor component satisfies
\[
 (\nabla_\mu A)(e_\alpha,e_\alpha)\ge0
 \quad\text{at }x_0.
\]

Suppose next that the null direction is normal, so that \(A_{nn}=0\).
By \eqref{eq:block},
\[
 A_{\alpha n}=0.
\]
Using \eqref{eq:Rfirst},
\[
 \rho\rho_\mu
 =
 h_iA_{in}
 =
 h_nA_{nn}+\sum_\alpha h_\alpha A_{\alpha n}
 =
 0.
\]
Since \(\rho>0\),
\[
 \rho_\mu=0.
\]
Consequently, the normal derivative of the radial Orlicz factor \(G(\rho)\) vanishes at \(x_0\).

Differentiating
\[
 \Psi(A)=\Hh
\]
in the normal direction yields
\[
 \Psi^{ij}A_{ij;n}=\Hh_\mu.
\]
By \eqref{eq:bdryA},
\[
 \sum_\alpha\Psi^{\alpha\alpha}A_{\alpha\alpha;n}
 =
 -\cot\theta
 \sum_\alpha\Psi^{\alpha\alpha}A_{\alpha\alpha}.
\]
Since \(A_{nn}=0\) and \(\Psi\) is one-homogeneous,
\[
 \sum_\alpha\Psi^{\alpha\alpha}A_{\alpha\alpha}
 =
 \Psi(A)=\Hh.
\]
Therefore
\[
 \Psi^{nn}A_{nn;n}
 =
 \Hh_\mu+\cot\theta\,\Hh.
\]

At the present point, \(\rho_\mu=0\), so
\[
 \frac{\Hh_\mu}{\Hh}
 =
 \frac1k\frac{f_\mu}{f}
 +
 \frac{p-1}{k}\frac{h_\mu}{h}
 =
 \frac1k\frac{f_\mu}{f}
 +
 \frac{p-1}{k}\cot\theta.
\]
Because
\[
 g=f^{-1/(p+k-1)},
\]
we have
\[
 \frac{f_\mu}{f}
 =
 -(p+k-1)\frac{g_\mu}{g}.
\]
The boundary inequality
\[
 g_\mu\le\cot\theta\,g
\]
therefore gives
\[
 \frac{f_\mu}{f}
 \ge
 -(p+k-1)\cot\theta.
\]
Hence
\[
 \frac{\Hh_\mu}{\Hh}\ge-\cot\theta,
\]
and so
\[
 \Hh_\mu+\cot\theta\,\Hh\ge0.
\]
Since \(\Psi^{nn}>0\),
\[
 A_{nn;n}\ge0.
\]

We have therefore shown that, if a boundary kernel exists, one can choose a unit null vector
\(v\) that is purely tangential or purely normal and satisfies
\[
 (\nabla_\mu A)(v,v)\ge0.
\]
We now make the opposite, strict inequality precise without differentiating an eigenvalue.
Put
\[
 w(x)=\lambda_{\min}(A(x)).
\]
By the interior full-rank conclusion, \(w>0\) in \(\Om\), whereas \(w\ge0\) on \(\overline\Om\).
Theorem~3.4 of \cite{BIS2020}, with \(m=1\), gives the following scalar inequality:
if \(\psi\) is any smooth lower support of \(w\) at an interior point, then, at the contact point,
\[
 \Psi^{ij}\psi_{;ij}\le C\bigl(\psi+|\nabla\psi|\bigr).
\]
Here a lower support means that \(\psi\le w\) locally with equality at that point.

For this boundary argument the constants need only be controlled for the fixed smooth solution.
Indeed, compactness and positivity give constants \(L,s_0>0\) such that
\[
 0\le A\le L\Id,
 \qquad
 \sigma_k(A)\ge s_0
 \quad\text{on }\overline\Om.
\]
Order the eigenvalues as \(\lambda_1\ge\cdots\ge\lambda_n\ge0\). Since
\[
 \sigma_k(A)\le\binom nk L^{k-1}\lambda_k,
 \qquad
 \lambda_k\ge\frac{s_0}{\binom nk L^{k-1}}>0,
\]
for every \(i\) the quantity \(\sigma_{k-1}(\lambda\mid i)\) has a positive lower bound.
The identity
\[
 \Psi^{ii}
 =\frac1k\sigma_k(A)^{1/k-1}\sigma_{k-1}(\lambda\mid i)
\]
therefore gives constants \(0<\lambda_*\le\Lambda_*<\infty\) with
\[
 \lambda_*\Id\le(\Psi^{ij})\le\Lambda_*\Id
 \quad\text{on }\overline\Om.
\]
This uses no positive lower bound for \(w\).
Moreover the matrices reached by this solution form a compact subset of the spectral domain \(\Gamma_k\),
and the fixed fields \(h,g,\Hh\) and their required derivatives are smooth and bounded up to the boundary,
with positive lower bounds for \(h,g,\Hh\).
Thus the structural calculation above and the local differential calculation in
\cite[Theorem~3.4]{BIS2020} have bounded coefficients and one constant \(C\)
throughout an interior collar of \(\partial\Om\).

Suppose now that \(w(x_0)=0\) at \(x_0\in\partial\Om\).
Take a sufficiently small interior tangent geodesic ball \(B_R(y)\subset\Om\), with
\(x_0\in\partial B_R(y)\), whose closure meets \(\partial\Om\) only at \(x_0\).
Take \(R<\pi/2\) and small enough that the ball lies in the collar just described.
Write \(r=d(y,\cdot)\), and on \(R/2<r<R\) set
\[
 b=e^{-q r^2}-e^{-q R^2},
 \qquad
 \mathcal L=\Psi^{ij}\nabla_{ij}.
\]
The spherical distance identity gives
\[
 \nabla^2b
 =e^{-q r^2}\left[
 4q^2r^2\,\dd r\otimes\dd r
 -2q\left\{\dd r\otimes\dd r
  +r\cot r\,(\sig-\dd r\otimes\dd r)\right\}\right].
\]
Since \(0<r\cot r\le1\) here, the ellipticity bounds imply
\[
 \mathcal Lb-C|\nabla b|-Cb
 \ge e^{-q r^2}
 \left[q^2\lambda_*R^2-2q(n\Lambda_*+CR)-C\right]>0
\]
when \(q\) is chosen sufficiently large.
Choose \(\varepsilon>0\) small enough that \(w\ge\varepsilon b\) on \(r=R/2\);
this also holds on \(r=R\), where \(b=0\).
If the minimum \(m_0\) of \(w-\varepsilon b\) on the closed annulus were negative,
it would be attained in its interior, where \(\psi=\varepsilon b+m_0\) is a smooth lower support of \(w\).
The lower-support inequality would give
\[
 \varepsilon\mathcal Lb
 \le C\left(\varepsilon b+m_0+\varepsilon|\nabla b|\right)
 <\varepsilon C\left(b+|\nabla b|\right),
\]
contradicting the strict inequality for \(b\). Hence \(w\ge\varepsilon b\) on this annulus.

Along the inward normal geodesic \(\gamma(s)=\exp_{x_0}(-s\mu)\), one has \(r=R-s\), and consequently
\[
 w(\gamma(s))\ge\varepsilon b(\gamma(s))\ge c s
 \quad\text{for all sufficiently small }s>0,
\]
for some \(c>0\).
Choose the pure tangential or pure normal unit null vector \(v\) found above,
and parallel transport it along \(\gamma\), writing the resulting field as \(v(s)\).
The smooth scalar component
\[
 a(s)=A_{\gamma(s)}(v(s),v(s))
\]
satisfies \(a(0)=0\) and \(a(s)\ge w(\gamma(s))\ge cs\).
Therefore
\[
 (\nabla_\mu A)(v,v)=-a'(0)\le-c<0,
\]
a contradiction. Thus no boundary null direction can occur.
Only the smooth tensor component \(a\) has been differentiated; no differentiability of \(w\) is assumed.

\begin{theorem}[Full rank]\label{thm:fullrank}
Let \(h>0\) be a smooth solution of \eqref{eq:intro-main}--\eqref{eq:intro-robin} with \(A[h]\ge0\).
Assume \eqref{eq:gconv}--\eqref{eq:gbdry}.
Then
\[
 A[h]>0
 \quad\text{on }\overline\Om.
\]
For this conclusion, no sign condition on \(Q_\varphi'\) is needed.
\end{theorem}

\subsection{Higher Regularity and Compactness}

By \cref{prop:c01,prop:c2}, there exist constants \(c_0,C_0,L>0\) such that
\[
 c_0\le h\le C_0,
 \qquad
 |\nabla h|\le C_0,
 \qquad
 0\le A[h]\le L\Id.
\]
Moreover,
\[
 \sigma_k(A[h])\ge s_0>0.
\]

\begin{lemma}\label[lemma]{lem:elliptic}
The operator \(\sigma_k^{1/k}\) is uniformly elliptic on the matrix set reached by the homotopy solutions.
\end{lemma}

\begin{proof}
Let
\[
 \lambda_1\ge\cdots\ge\lambda_n\ge0.
\]
Since
\[
 s_0\le\sigma_k(\lambda)\le\binom nkL^{k-1}\lambda_k,
\]
we obtain \(\lambda_k\ge c>0\).
Thus every derivative
\[
 \frac{\partial\sigma_k^{1/k}}{\partial\lambda_i}
\]
has a positive lower bound, while the upper bound follows from \(\lambda_i\le L\).
\end{proof}

To apply global oblique regularity, use the concave root of the equation.
In smooth local coordinates on the fixed spherical cap, set
\[
 \begin{aligned}
 \mathcal A_{ij}(x,z,v,R)
 &=R_{ij}-\Gamma_{ij}^{a}(x)v_a+z\sig_{ij}(x),\\
 F_t(x,z,v,R)
 &=\sigma_k^{1/k}(\mathcal A)-\widehat f_t(x,z,v)^{1/k},
 \end{aligned}
\]
where the eigenvalues in \(\sigma_k\) are taken with respect to \(\sig\).
Thus \(F_t(x,h,Dh,D^2h)=0\), and \(F_t\) is concave in \(R\).
By \cref{prop:c01,prop:c2}, all lower-order jets lie in a fixed compact region
with \(z\ge c_0>0\) and bounded \(z\) and \(v\); the corresponding radial
variable \((z^2+|v|_{\sig}^2)^{1/2}\) therefore lies in a compact subinterval
of \((0,\infty)\).  On a slightly larger such region,
\eqref{eq:gt}--\eqref{eq:phit} and the positivity and smoothness of the fixed
data give uniform positive lower bounds for \(\widehat f_t\) and uniform
bounds for all its derivatives of any fixed finite order, independently of
\(t\in[0,1]\).

The attained matrices belong to the compact set
\[
 \mathcal K
 =\{A:0\le A\le L\Id,\ \sigma_k(A)\ge s_0\}
 \Subset\Gamma_k.
\]
Here \(\Gamma_k=\{A:\sigma_j(A)>0,\ 1\le j\le k\}\) is the matrix
admissibility cone.
By \cref{lem:elliptic}, the concave operator \(\sigma_k^{1/k}\) is uniformly
elliptic on a fixed neighborhood of \(\mathcal K\) within \(\Gamma_k\).
Consequently \(F_t\), localized to a fixed neighborhood of the attained
jets, has uniformly controlled smooth coefficients and uniform ellipticity
constants.  The cap and its boundary are smooth.  In boundary charts the
Robin operator is
\[
 B(x,z,v)=\mu^i(x)v_i-\cot\theta\,z,
\]
which is smooth and uniformly oblique; to use the inward-normal convention
of \cite{LiebermanTrudinger1986}, take the boundary operator \(-B\).
The bounded-\(C^2\) global Hessian H\"older estimate of
\cite[Section~6, Theorem~6.2, together with the supplementary remarks on
p.~544]{LiebermanTrudinger1986}, applied in this localized jet region and
in smooth boundary-flattening charts, therefore gives
\[
 \|h\|_{C^{2,\alpha}(\overline\Om)}\le C
\]
for some fixed \(\alpha\in(0,1)\), with \(C\) independent of \(t\).
Here the supplementary remarks allow Hessian concavity in place of
condition (6.2) and remove condition (6.3) of Section~6.  All coefficient,
ellipticity, obliqueness, and chart bounds used in this estimate are uniform
in \(t\); no positive lower bound for the smallest eigenvalue of \(A[h]\)
is used.
The linearized coefficients are now \(C^\alpha\), so the higher-order oblique
Schauder estimates, with the same uniform coefficient control, yield
\begin{equation}\label{eq:C4}
 \|h\|_{C^{4,\alpha}(\overline\Om)}\le C,
\end{equation}
where \(C\) is independent of \(t\).

We now obtain a uniform strict-convexity margin.
Assume otherwise. Then there exist \(t_j\in[0,1]\), corresponding homotopy solutions \(h_j\), and points
\(x_j\in\overline\Om\) such that the smallest eigenvalue of \(A[h_j](x_j)\) tends to zero.
Passing to a subsequence,
\[
 t_j\to t_\infty.
\]
By \eqref{eq:C4} and compact embedding, for every \(0<\beta<\alpha\) we may assume
\[
 h_j\to h_\infty
 \quad\text{in }C^{4,\beta}(\overline\Om).
\]
The uniform \(C^0\)-estimate gives
\[
 h_\infty>0,
\]
while the \(C^2\)-convergence gives
\[
 A[h_\infty]\ge0.
\]
Passing to the limit in the equation and in the Robin condition shows that
\(h_\infty\) is a solution of the homotopy problem with parameter
\(t_\infty\).  Its positive right-hand side and the bounds above put
\(A[h_\infty]\) in the same compact subset of \(\Gamma_k\), even though
only \(A[h_\infty]\ge0\) is known at this stage.  The preceding uniformly
elliptic oblique Schauder bootstrap therefore gives
\(h_\infty\in C^\infty(\overline\Om)\) before any strict-convexity conclusion
is invoked.  The structural conditions on \(g_t\) are preserved along the
homotopy, so \cref{thm:fullrank} applies:

\[
 A[h_\infty]>0
 \quad\text{on }\overline\Om.
\]
Since \(\overline\Om\) is compact, the smallest eigenvalue of \(A[h_\infty]\) has a positive minimum.
The \(C^2\)-convergence then implies that \(A[h_j]\) is uniformly positive for all sufficiently large \(j\),
contradicting the choice of \(h_j\) and \(x_j\).
Hence there exists \(\kappa>0\), independent of \(t\), such that
\begin{equation}\label{eq:kappa}
 A[h]\ge\kappa\sig
\end{equation}
for every homotopy solution.

\section{Proof of the Main Theorem}

\subsection{Proof of Theorem 1.1}

Fix \(0<\beta<\alpha\) and let
\[
 X^+=C_{\mathrm{even}}^{4,\beta}(\overline\Om),
 \qquad
 Y^+=C_{\mathrm{even}}^{2,\beta}(\overline\Om)
       \times C_{\mathrm{even}}^{3,\beta}(\partial\Om).
\]
Choose \(M\) so large that the uniform \(C^{4,\alpha}\)-estimate gives
\(\|h\|_{C^{4,\beta}}<M\) for every homotopy solution. Define
\begin{equation}\label{eq:U}
 \U
 =
 \left\{
 h\in X^+:
 \|h\|_{C^{4,\beta}}<M,\;
 h>\frac{c_0}{2},\;
 A[h]>\frac{\kappa}{2}\sig
 \right\}.
\end{equation}
Then \(\U\) is a bounded open subset of \(X^+\). All boundaries of \(\U\) below are taken
relative to \(X^+\).

Define
\[
 \mathcal P_t(h)
 =
 \sigma_k(A[h])
 -
 f_t h^{p-1}\frac{\rho^{k+1}}{\varphi_t(\rho)},
 \qquad
 \mathcal B(h)=\nabla_\mu h-\cot\theta\,h.
\]
Thus \((\mathcal P_t,\mathcal B):\U\to Y^+\). On a neighborhood of \(\overline\U\),
the interior operator is uniformly elliptic and the boundary operator is uniformly oblique,
with coefficient bounds uniform in \(t\).

We make explicit the compact construction used to define the degree on \(X^+\).
The gnomonic chart \(\xi\mapsto\xi'/\xi_{n+1}\) identifies \(\overline\Om\) smoothly
with the closed Euclidean ball of radius \(\tan\theta\), and the horizontal reflection becomes
\(y\mapsto-y\). We suppress the coordinate pullbacks. In these coordinates, write
\[
 \begin{aligned}
 X&=C^{4,\beta}(\overline\Om),\\
 Y&=C^{2,\beta}(\overline\Om)\times C^{3,\beta}(\partial\Om),\\
 Z&=C^\beta(\overline\Om)\times C^{1,\beta}(\partial\Om)
                            \times C^{1,\beta}(\partial\Om).
 \end{aligned}
\]
Let \(\Delta_E\) and \(\nu\) denote the Euclidean Laplacian and outer unit normal in this chart.
The auxiliary isomorphism \(\mathcal D:Y\to Z\) is
\[
 \mathcal D(r,s)
 =\bigl(\Delta_E r,\;\partial_\nu r+r,\;
                    \Delta_{\partial\Om}s-s\bigr),
\]
where \(\Delta_{\partial\Om}\) is the Euclidean tangential Laplacian.
Let \(\mathcal O\subset X\) be the ambient open set defined by the inequalities in
\eqref{eq:U} without the symmetry condition.
To match the global jet-domain formulation in \cite{LLN2015}, extend the lower-order term
smoothly to all scalar and gradient jets using a reflection-invariant cutoff that equals one
near the lower jets attained on \(\overline{\mathcal O}\) and is supported where \(z>c_0/4\).
The extension is uniform in \(t\) and changes neither the operators nor their principal
coefficients near \(\overline\U\).
Apply the construction in \cite[Section~3, Theorem~2]{LLN2015} on \(\mathcal O\).
For a common sufficiently large \(N\), its frozen fourth-order operator
\(\mathcal L_{h,N}:X\to Z\) is an isomorphism. Since the principal coefficients of
\(\mathcal P_t\) and \(\mathcal B\) are independent of \(t\), so is \(\mathcal L_{h,N}\).
Set
\[
 Q_h=\mathcal L_{h,N}^{-1}\mathcal D:Y\longrightarrow X,
 \qquad
 \Phi_t(h)=Q_h\bigl(\mathcal P_t(h),\mathcal B(h)\bigr)=h+K_t(h).
\]
The decomposition and estimates in \cite[Section~3]{LLN2015} show that \(K_t\) is a compact
homotopy. All these maps extend continuously to \(\overline\U\).

Let \(J\) denote pullback by \(y\mapsto-y\), acting componentwise on the data spaces.
The Euclidean auxiliary operators commute with \(J\), and covariance of the frozen principal
coefficients gives
\[
 \mathcal D J=J\mathcal D,\qquad
 \mathcal L_{Jh,N}J=J\mathcal L_{h,N},\qquad
 Q_{Jh}J=JQ_h.
\]
The evenness of the prescribed data also gives
\[
 \Phi_t(Jh)=J\Phi_t(h),\qquad K_t(Jh)=JK_t(h).
\]
Consequently, \(Q_h:Y^+\to X^+\) is an isomorphism for every \(h\in\overline\U\), and
\(K_t\) restricts to a compact homotopy on the closed subspace \(X^+\).
We therefore define
\[
 \deg_+\bigl((\mathcal P_t,\mathcal B),\U,0\bigr)
 :=\deg_{\mathrm{LS},X^+}(\Phi_t,\U,0).
\]

The uniform estimates show that every zero of the homotopy pair in \(\overline\U\) satisfies
\[
 h\ge c_0>\frac{c_0}{2},\qquad
 A[h]\ge\kappa\sig>\frac{\kappa}{2}\sig,\qquad
 \|h\|_{C^{4,\beta}}<M.
\]
Since each \(Q_h\) is invertible, this excludes zeros of \(\Phi_t\) on the relative boundary
\(\partial\U\), uniformly for \(t\in[0,1]\). Homotopy invariance gives
\begin{equation}\label{eq:degree-homotopy}
 \deg_+((\mathcal P_1,\mathcal B),\U,0)
 =
 \deg_+((\mathcal P_0,\mathcal B),\U,0).
\end{equation}

At \(t=0\),
\[
 \varphi_0(r)=r^{k+1},\qquad
 f_0=\binom nk\ell^{1-p},\qquad h_0=\ell.
\]
By \eqref{eq:ell}, \(h_0\) is a solution and hence belongs to \(\U\).
The uniqueness theorem \cite[Theorem~5.4]{MWWLp2025}, with \(l=0\), implies that \(\ell\)
is the only zero of \((\mathcal P_0,\mathcal B)\) in \(\U\), since every such zero is
positive and strictly convex and \(1<p<k+1\).
The linearization is
\[
 D\mathcal P_0(\ell)v
 =\binom{n-1}{k-1}
   \left(\Delta v+nv-\frac{n(p-1)}{k\ell}v\right),
 \qquad
 D\mathcal B(\ell)v=\mathcal B(v).
\]
By \cite[Lemma~5.6]{MWWLp2025}, \(D\mathcal P_0(\ell)\) is an isomorphism from
\(\{v\in X^+:\mathcal B(v)=0\}\) onto \(C_{\mathrm{even}}^{2,\beta}(\overline\Om)\).
This also gives an isomorphism for the full pair
\[
 D(\mathcal P_0,\mathcal B)(\ell):X^+\longrightarrow Y^+.
\]
Indeed, for any even boundary datum \(s\in C^{3,\beta}(\partial\Om)\), solve the standard
Neumann problem
\[
 \Delta z_s-z_s=0\quad\text{in }\Om,\qquad
 \nabla_\mu z_s=s/\ell\quad\text{on }\partial\Om.
\]
Its unique solution belongs to \(C^{4,\beta}(\overline\Om)\) and is even by uniqueness.
Then \(v_s=\ell z_s\in X^+\) satisfies \(\mathcal B(v_s)=s\), by \eqref{eq:ell}.
The homogeneous-Robin isomorphism now uniquely solves
\[
 D\mathcal P_0(\ell)w=r-D\mathcal P_0(\ell)v_s,\qquad \mathcal B(w)=0
\]
for every even \(r\in C^{2,\beta}(\overline\Om)\). Thus \(v_s+w\) solves the full pair;
injectivity follows from the same homogeneous-Robin isomorphism.

Since \((\mathcal P_0(\ell),\mathcal B(\ell))=0\), differentiating the compact degree map gives
\[
 D\Phi_0(\ell)=Q_\ell D(\mathcal P_0,\mathcal B)(\ell):X^+\longrightarrow X^+,
\]
an isomorphism. Choose a sufficiently small ball \(V\) about \(\ell\), with
\(\overline V\subset\U\). The local linearization property of Leray--Schauder degree gives
\(\deg_{\mathrm{LS},X^+}(\Phi_0,V,0)=\pm1\).
Uniqueness of the zero and excision therefore give
\[
 \deg_+((\mathcal P_0,\mathcal B),\U,0)
 =\deg_{\mathrm{LS},X^+}(\Phi_0,V,0)=\pm1.
\]
By \eqref{eq:degree-homotopy}, the degree at \(t=1\) is nonzero. Its existence property
provides \(h\in\U\) with
\[
 \mathcal P_1(h)=0,\qquad \mathcal B(h)=0.
\]
Thus \(h\) solves the target equation and satisfies \(h>0\) and \(A[h]>0\).
Finally, the higher regularity argument gives \(h\in C^\infty(\overline\Om)\).
This proves \cref{thm:main}.

\subsection{Examples and Remarks}

\begin{example}
Let
\[
 \varphi(r)=r^q.
\]
Then
\[
 H(r)=r^{q-p},
 \qquad
 Q_\varphi(r)=q.
\]
The assumptions on \(\varphi\) in \cref{thm:main} reduce exactly to
\[
 p<q\le k+1.
\]
Thus the power model is contained in the present Orlicz theorem.
\end{example}

\begin{example}
Let
\[
 \varphi(r)=r^q(1+r)^\beta,
\]
where
\[
 p<q,\qquad \beta\ge0,\qquad q+\beta\le k+1.
\]
Then
\[
 Q_\varphi(r)
 =
 q+\frac{\beta r}{1+r},
\]
and therefore
\[
 p<Q_\varphi(r)\le k+1.
\]
Moreover,
\[
 H(r)=r^{q-p}(1+r)^\beta
\]
is strictly increasing and proper.
Hence \cref{thm:main} applies.
\end{example}

\begin{remark}
The assumptions on \(\varphi\) have separate roles.
The proper monotonicity of \(H=\varphi/r^p\) is a scale condition.
The inequality \(Q_\varphi\le k+1\) is used only in the weighted gradient estimate.
The full-rank theorem itself does not require a sign condition on \(Q_\varphi'\).
\end{remark}

\begin{remark}
The evenness assumption enters only through the currently available capillary non-collapsing theorem.
All local estimates and the full-rank argument are non-symmetric.
\end{remark}

\begin{remark}
The present argument proves existence, not uniqueness, for a general Orlicz function.
The nonhomogeneity of \(\varphi\) makes the uniqueness problem substantially different from the power case.
This is analogous to the broader dual Orlicz setting, where uniqueness is typically much more delicate than existence;
see \cite{ZXY2018}.
\end{remark}

\medskip
\noindent\textbf{Acknowledgments:}
Hua Zhu was supported by the National Natural Science Foundation of
China (Grant Nos. 12661041 and 12501273 ) and the Research Foundation
of Southwest University of Science and Technology
(Grant Nos. 25zx7154 and 25zx7153). Qiqi Zhang was supported by the National Natural Science Foundation of
China (Grant Nos. 12301098 and 12661014 ).
\medskip

\section*{Declaration on generative AI assistance}

During the development and preparation of this work, the authors used a generative AI tool
for preliminary, non-authoritative exploratory assistance, including organizational discussion
and the consideration of possible approaches. The AI-generated outputs were not treated as
mathematical sources, and no argument was included without independent verification and
substantial revision by the authors. All theorem statements, proofs, computations, references,
and final arguments were independently checked, revised, and finalized by the authors, who
take full responsibility for the correctness, originality, and integrity of the paper.


\begin{thebibliography}{99}

\bibitem{BIS2020}
P. Bryan, M.N. Ivaki, J. Scheuer,
Constant rank theorems for curvature problems via a viscosity approach.
\emph{Calc. Var. Partial Differential Equations} \textbf{62}, Art.~98 (2023).

\bibitem{CabezasHu2025}
C. Cabezas-Moreno, J. Hu,
The \(L_p\) dual Christoffel--Minkowski problem for \(1<p<q\le k+1\) with \(1\le k\le n\).
\emph{Calc. Var. Partial Differential Equations} \textbf{64}, Art.~229 (2025).

\bibitem{ChouWang2006}
K.S. Chou, X.J. Wang,
The \(L_p\)-Minkowski problem and the Minkowski problem in centroaffine geometry.
\emph{Adv. Math.} \textbf{205}, 33--83 (2006).

\bibitem{DingLi2025}
S. Ding, G. Li, 
A type of anisotropic flows and dual Orlicz Christoffel--Minkowski type equations.
\emph{J. Differential Equations} \textbf{445}, Art.~113586 (2025).

\bibitem{Gardner1995}
R.J. Gardner,
\emph{Geometric Tomography}.
Cambridge University Press, Cambridge (1995).

\bibitem{GardnerEtAl2014}
R.J. Gardner, D. Hug, W. Weil,
The Orlicz--Brunn--Minkowski theory: a general framework, additions, and inequalities.
\emph{J. Differential Geom.} \textbf{97}, 427--476 (2014).

\bibitem{GardnerEtAl2015}
R.J. Gardner, D. Hug, W. Weil, D. Ye,
The dual Orlicz--Brunn--Minkowski theory.
\emph{J. Math. Anal. Appl.} \textbf{430}, 810--829 (2015).

\bibitem{HaberlEtAl2010}
C. Haberl, E. Lutwak, D. Yang, G. Zhang,
The even Orlicz Minkowski problem.
\emph{Adv. Math.} \textbf{224}, 2485--2510 (2010).

\bibitem{HuMaShen2004}
C. Hu, X.-N. Ma, C. Shen,
On the Christoffel--Minkowski problem of Firey's \(p\)-sum.
\emph{Calc. Var. Partial Differential Equations} \textbf{21}, 137--155 (2004).

\bibitem{HY2026}
J. Hu, B. Yang,
Capillary John ellipsoid theorem with applications to capillary curvature problems, arXiv:2603.27252v1(2026).

\bibitem{HuIvaki2026}
Y. Hu, M.N. Ivaki,
Capillary \(L_p\)-curvature problem.
\emph{Adv. Calc. Var.} \textbf{19}(4), 699--705 (2026).


\bibitem{HuIvakiScheuer2025}
Y. Hu, M.N. Ivaki, J. Scheuer,
Capillary Christoffel--Minkowski problem.
{arXiv:2504.09320v1} (2025).

\bibitem{HuangHe2012}
Q. Huang, B. He,
On the Orlicz Minkowski problem for polytopes.
\emph{Discrete Comput. Geom.} \textbf{48}, 281--297 (2012).

\bibitem{HLYZ2016}
Y. Huang, E. Lutwak, D. Yang, G. Zhang,
Geometric measures in the dual Brunn--Minkowski theory and their associated Minkowski problems.
\emph{Acta Math.} \textbf{216}, 325--388 (2016).

\bibitem{Li2014}
A.J. Li,
The generalization of Minkowski problems for polytopes.
\emph{Geom. Dedicata} \textbf{168}, 245--264 (2014).

\bibitem{LiLiu2026}
G. Li, M. Liu,
The capillary \(L_p\) dual Christoffel--Minkowski problem for \(1<p<q\le k+1\) with \(1\le k\le n\), arXiv:2609.37258v1(2026).

\bibitem{LLN2015}
Y.Y. Li, J. Liu, L. Nguyen,
A degree theory for second order nonlinear elliptic operators with nonlinear oblique boundary conditions.
\emph{J. Fixed Point Theory Appl.} \textbf{19}, 853--876 (2017).


\bibitem{LiebermanTrudinger1986}
G.M. Lieberman, N.S. Trudinger,
Nonlinear oblique boundary value problems for nonlinear elliptic equations.
\emph{Trans. Amer. Math. Soc.} \textbf{295}, 509--546 (1986).

\bibitem{Lutwak1975}
E. Lutwak,
Dual mixed volumes.
\emph{Pacific J. Math.} \textbf{58}, 531--538 (1975).

\bibitem{Lutwak1988}
E. Lutwak,
Intersection bodies and dual mixed volumes.
\emph{Adv. Math.} \textbf{71}, 232--261 (1988).

\bibitem{Lutwak1993}
E. Lutwak, 
The Brunn--Minkowski--Firey theory I: mixed volumes and the Minkowski problem.
\emph{J. Differential Geom.} \textbf{38}, 131--150 (1993).

\bibitem{LYZ2004}
E. Lutwak, D. Yang, G. Zhang,
On the \(L_p\)-Minkowski problem.
\emph{Trans. Amer. Math. Soc.} \textbf{356}, 4359--4370 (2004).

\bibitem{MWW2025}
X. Mei, G. Wang, L. Weng,
The capillary Christoffel--Minkowski problem.
\emph{Calc. Var. Partial Differential Equations} \textbf{65}, Art.~186 (2026).

\bibitem{MWWLp2025}
X. Mei, G. Wang, L. Weng,
Prescribed \(L_p\) curvature problem for convex capillary hypersurface, arXiv:2512.16686v3(2026).

\bibitem{MWX2024}
X. Mei, G. Wang, L. Weng, C. Xia,
Alexandrov--Fenchel inequalities for convex hypersurfaces in the half-space with capillary boundary II.
\emph{Math. Z.} \textbf{310}, Art.~71 (2025).


\bibitem{Schneider2014}
R. Schneider,
\emph{Convex Bodies: The Brunn--Minkowski Theory}, 2nd expanded ed.
Cambridge University Press, Cambridge (2014).

\bibitem{WWX2022}
G. Wang, L. Weng, C. Xia,
Alexandrov--Fenchel inequalities for convex hypersurfaces in the half-space with capillary boundary.
\emph{Math. Ann.} \textbf{388}, 2121--2154 (2024).

\bibitem{XiJinLeng2014}
D. Xi, H. Jin, G. Leng, The Orlicz Brunn--Minkowski inequality.
\emph{Adv. Math.} \textbf{260}, 350--374 (2014).

\bibitem{Zhao2015}
C.J. Zhao, Orlicz dual mixed volumes.
\emph{Results Math.} \textbf{68}, 93--104 (2015).

\bibitem{Zhao2017}
Y. Zhao, The dual Minkowski problem for negative indices.
\emph{Calc. Var. Partial Differential Equations} \textbf{56}, Art.~18 (2017).

\bibitem{ZXY2018}
B. Zhu, S. Xing, D. Ye,
The dual Orlicz--Minkowski problem.
\emph{J. Geom. Anal.} \textbf{28}, 3829--3855 (2018).

\bibitem{ZhuZhouXu2014}
B. Zhu, J. Zhou, W. Xu,
Dual Orlicz--Brunn--Minkowski theory.
\emph{Adv. Math.} \textbf{264}, 700--725 (2014).

\bibitem{Zhu2026}
H. Zhu,
\emph{The Capillary Orlicz--Christoffel--Minkowski Problem},
arXiv:2611.18299v2(2026).


\end{thebibliography}
\end{document}